\documentclass[11pt]{article}

\usepackage{graphicx}
\usepackage[centertags]{amsmath}
\usepackage{amsfonts}
\usepackage{amsmath}
\usepackage{amssymb}
\usepackage{latexsym}
\usepackage{amsthm}
\usepackage{mathrsfs}
\usepackage{newlfont}
\usepackage[latin1]{inputenc}

\numberwithin{equation}{section}

\title{Principal dynamical  components}

\author{Manuel D. de la Iglesia and Esteban G. Tabak \thanks{Courant Institute of Mathematical Sciences, 251 Mercer St.
New York, NY 10012, USA, {\tt mdi29@cims.nyu.edu, tabak@cims.nyu.edu}.}}

\begin{document}

\maketitle

\begin{abstract}
A new procedure is proposed for the dimensional reduction of time series. Similarly to principal components, the procedure seeks a low-dimensional manifold that minimizes information loss. Unlike principal components, however, the new procedure involves dynamical considerations, through the proposal of a predictive dynamical model in the reduced manifold. Hence the minimization of the uncertainty is not only over the choice of a reduced manifold, as in principal components, but also over the parameters of the dynamical model. Further generalizations are provided to non-autonomous and non-Markovian scenarios, which are then applied to historical sea-surface temperature data.
\end{abstract}



\section{Introduction}

Complex systems typically involve a large number of degrees of freedom. Thus to elucidate the fundamental mechanisms underlying one such system's behavior, one may consider its projection onto smaller-dimensional manifolds, selected so as to capture as much of the dynamics as possible. A tool frequently used for this purpose is principal components \cite{Jit}, whereby a linear subspace of prescribed dimensionality of the phase-space of observations is sought, so as to maximize the amount of the variability that is preserved when the data are projected onto it.

Given a dataset $z_j$, $j \in [1, \ldots, N]$, where each observation $z_j$ consists of $n$ real numbers, its first $m$ ($m \le n$) principal components are given by $x_j= Q_x' (z_j-\bar{z})$, where $\bar{z}$ is the mean value of $z$, and $Q_x$ is an $n \times m$ matrix with orthonormal columns, chosen so that $\sum_{j=1}^N \left\| \left(z_j-\bar{z}\right)-Q_x x_j\right\|^2$ is as small as possible. From a statistical perspective, among all  $m$-dimensional  subspaces,  $x$ is the one whose knowledge minimizes the uncertainty of $z$. The matrix $Q_x$ consists of the first $m$ columns of $U$ in the singular value decomposition
$$ Z = U S V' , $$
where the elements $Z_i^j$ of the matrix $Z \in R^{n\times N}$ contain the $i$th component of the $j$th observation minus its mean value $\bar{z}^i$ over all observations, $U \in R^{n\times n}$ and $V \in R^{N \times N}$ are orthogonal matrices, and $S \in R^{n \times N}$ is the diagonal matrix of singular values of $Z$, the eigenvalues of the empirical covariance matrix $C = Z Z'$ sorted in decreasing order.

In the probabilistic scenario underlying this procedure, the $z_j$'s are independent samples of a Gaussian distribution ${\cal N}(\mu,\Sigma)$, $\bar{z}$ is an estimate for its mean $\mu$, and the principal components estimate the principal axes of the covariance matrix $\Sigma$, sorted in decreasing order by the fraction of the total variance that they explain. Yet principal components are often sought for data that do not quite fit this scenario. Of particular concern to us here is the situation where the $z_j$'s form a  time series, representing snapshots of the vector $z$ at  equidistant times $t_j$. In this context, the dimensional reduction by principal components, oriented toward data compression, lacks any concept of \emph{dynamics}: the various snapshots $z_j$ are treated as independent observations, which renders immaterial even the order in which they are sorted. If there is an underlying dynamics, this is neither unveiled nor exploited by the analysis.

An example is provided by the Empirical Orthogonal Functions (EOFs) \cite{EOF} --the name given to principal components in climate studies--, which take a time series of atmospheric or oceanic data, subtract its time average or ``climatology'', and find those modes that explain the largest share of its variability. These modes may then be assigned suggestive names such as ``El Ni\~no'' or ``The north-Atlantic oscillation'' and given a dynamical interpretation. Yet no dynamics ever entered into their calculation: just the static variability of the data, treated as a series of independent, unsorted observations.

In this paper, we develop an alternative methodology, highly reminiscent of the principal-component framework, but with a dynamical core. We seek, as in principal components, a hierarchy of manifolds, that we name ``principal dynamical components''. Attached to these manifolds is a model of predictive dynamics. The cost function to minimize has, as in principal components, the variability in the unrepresented variables, but also the fraction of the variability in the preserved variables that is not explained by the dynamics. Thus the dynamical components are characterized not by capturing most of the system's variability, but by explaining dynamically its largest possible share.
Hence this methodology can be thought of as a blend of autoregression analysis \cite{Wei}, which is used as a reduced dynamical model, and principal components, though the criterium for selecting a reduced manifold differs from the latter's.

We first present the new methodology, as a natural extension of  principal components, in a linear, autonomous framework, with a dynamic manifold given by $x= Q_x' z$, where $Q_x$ is a fixed  $n \times m$ orthogonal matrix, and the dynamics by $x_{j+1} = A x_j$, where $A$ is another fixed, $m\times m$ matrix. The definition of principal dynamical components results in a minimization problem over both $Q_x$ and $A$. Section \ref{Linear} presents this problem and provides an efficient methodology to solve it. Yet many real problems are not autonomous: climate dynamics, for instance, is season-dependent. In Section \ref{NonAutonomous} we extend the methodology to non-autonomous situations and, more generally, to accommodate for the presence of exogenous variables and external controls, that appear in many engineering applications. Here $Q_x$ and $A$ depend on time and on those external variables. Section \ref{NonM} extends the procedure further to handle non-Markovian processes, where the dynamics involves more than the immediate past. We illustrate the procedure throughout with synthetic data and, in Section \ref{Realex}, we concern ourselves with a real application to time series of sea-surface temperature over the ocean. Section \ref{Nonlinear} gives a probabilistic interpretation of  the principal  dynamical component procedure, which provides a conceptual extension to general nonlinear, non-Gaussian settings. The development of effective algorithms for the numerical implementation of this broad generalization will be described elsewhere.

\section{The linear, autonomous framework}
\label{Linear}

The probabilistic set-up for principal component analysis consists of independent observations drawn from a Gaussian distribution. The natural extension to time series has  a time series $z_j$ , $j \in [1,\ldots, N]$ drawn from the linear Markovian dynamics
$$ z_{j+1} =  \mathcal{N}(A^z z_j, \Sigma^z) . $$
Here the matrix $A^z$ models autocorrelation, and $\mathcal{N}$ represents a Gaussian process with mean $A^z z_j$ and covariance matrix $\Sigma^z$. Neither $A^z$ nor $\Sigma^z$ are known to us; instead, we seek an $m$-dimensional manifold $x= Q_x' z$ and reduced dynamics
$$ x_{j+1} = A x_j + \xi_j, $$
where $\xi_j$ is the prediction error, such that the \emph{predictive} uncertainty or cost
$$ c = \sum_{j=1}^{N-1} \left\|z_{j+1}-Q_x x_{j+1}\right\|^2  = \sum_{j=1}^{N-1} \left\|z_{j+1}-Q_x A Q_x' z_j \right\|^2 $$
is minimal. This is the conceptual basis of what we shall denote \emph{linear autonomous principal dynamical component} analysis.

It is convenient to introduce some further notation: $y$ for the orthogonal complement of $x$, so that
$$ z = [Q_x Q_y] \left( \begin{array}{c}
		x \\
		y
	\end{array}
\right) \, ,
$$
where $Q=[Q_x Q_y]$ is an orthogonal matrix, and $\tilde{x}$ for the conditional expectation of $x$:
$$ \tilde{x}_{j+1} = A x_j . $$
Since the dynamics of $y$ is not explained by the model, we have $\tilde{y}_{j+1}=0$.

\subsection{Two-dimensional case}

The simplest scenario, appropriate for a first view of the proposed algorithm, has the observations $z_j$ in a two-dimensional space, $n=2$, and seeks a reduced manifold $x$ of dimension $m=1$. We introduce the following notation:
\begin{equation*}
z = \left( \begin{array}{c}
		A \\
		P
	\end{array}
\right) \, ,
\end{equation*}
where, mimicking an application to climate dynamics, $A$ stands for Atlantic and $P$ for Pacific spatially-averaged sea-surface temperatures,
\begin{equation*}
x = A \cos(\theta_*) + P \sin(\theta_*) \, ,
\end{equation*}
\begin{equation*}
y = -A \sin(\theta_*) + P \cos(\theta_*) \, ,
\end{equation*}
where the angle $\theta_*$ defines the direction of the dynamic component $x$ in $(A,P)$ space, and
\begin{equation*}
\tilde{x}_{j+1} = a\, x_j \, ,
\end{equation*}
where the stretching factor $a$ describes the deterministic component of the reduced dynamics.

The cost function adopts the form
\begin{eqnarray*}
c(\theta, a) &=& \sum_{j=1}^{N-1} \left\|  \left( \begin{array}{c}
		A_{j+1} - \tilde{A}_{j+1}\\
		P_{j+1} - \tilde{P}_{j+1}
	\end{array}
\right)   \right\|^2 =
\sum_{j=1}^{N-1} \left\|  \left( \begin{array}{c}
		x_{j+1} - \tilde{x}_{j+1}\\
		y_{j+1} - \tilde{y}_{j+1}
	\end{array}
\right)   \right\|^2 \\ &=&
\sum_{j=1}^{N-1} \left\|  \left( \begin{array}{c}
		x_{j+1} - a x_j\\
		y_{j+1}
	\end{array}
\right)   \right\|^2 = \sum_{j=1}^{N-1} \left(y_{j+1}\right)^2 + \left(x_{j+1}-a x_j\right)^2 .
\end{eqnarray*}
By contrast, the corresponding cost function for regular principal components in this 2-dimensional scenario is
$$ c_{pc}(\theta) =  \sum_{j=1}^{N} \left(y_j\right)^2 : $$
the amount of variability in the unrepresented variable $y$.

The minimization of $c$ can be solved iteratively. If at the beginning of a step we have coordinates $(x,y)$, then
\begin{equation*}
\frac{\partial c}{\partial a} = -2 \sum_{j=1}^{N-1} \left(x_{j+1}-a x_j\right) x_j .
\end{equation*}
Equating $\frac{\partial c}{\partial a}$  to zero yields the standard regression formula
$$ a = \frac{\sum_{j=1}^{N-1}x_j x_{j+1}}{\sum_{j=1}^{N-1} x_j^2} . $$
If now we update $x$ and $y$ through a further rotation
\begin{equation*}
x \leftarrow x \cos(\theta) + y \sin(\theta) \, ,
\end{equation*}
\begin{equation*}
y \leftarrow -x \sin(\theta) + y \cos(\theta) \, ,
\end{equation*}
we have
\begin{eqnarray*}
\frac{\partial c}{\partial \theta} = 2  a \sum_{j=1}^{N-1}
\Big[\left(a x_j y_j - \left(x_{j+1} y_j + x_j y_{j+1} \right) \right) \cos(2 \theta) + \\
\left(x_{j+1}x_j-y_{j+1}y_j + \frac{a}{2} \left(y_j^2-x_j^2\right) \right) \sin(2 \theta) \Big] .
\end{eqnarray*}
Rather than seeking a  closed expression for $\theta$ that would make this derivative vanish  --notice that $\theta$ is implicitly included in the definition of the $x$ and $y$'s--, it is preferable to descend the gradient
\begin{equation*}
\frac{\partial c}{\partial \theta}\Big|_{\theta=0} = 2 a \sum_{j=1}^{N-1}
\left[a x_j y_j - \left(x_{j+1} y_j + x_j y_{j+1} \right) \right]
\end{equation*}
or, more efficiently, to involve also the second derivative
\begin{equation*}
\frac{\partial^2 c}{\partial \theta^2}\Big|_{\theta=0} = 2 a \sum_{j=1}^{N-1}
\left[2 x_{j+1}x_j-2 y_{j+1}y_j + a \left(y_j^2-x_j^2\right) \right],
\end{equation*}
and compute the $\theta$ that minimizes the quadratic local approximation to $c$:
$$ \theta = \theta_q = -\frac{ \frac{\partial c}{\partial \theta}\Big|_{\theta=0} }{\frac{\partial^2 c}{\partial \theta^2}\Big|_{\theta=0} }.$$
A little extra care is required when applying the quadratic approximation far from the optimal $\theta$: if $\frac{\partial^2 c}{\partial \theta^2}\Big|_{\theta=0}  \le 0$, then we must do descent instead:
\begin{equation}
  \theta = -\epsilon_l \frac{\partial c}{\partial \theta}\Big|_{\theta=0},
  \label{descent}
\end{equation}
where $\epsilon_l > 0$ is a chosen learning rate. Also, if $\theta_q$ is too big, we must limit our step size:
\begin{equation*}
  |\theta| = \max(|\theta_q|, \epsilon),
\end{equation*}
where $\epsilon$ is the maximum allowable step in $\theta$. It is sensible to relate the values of the two $\epsilon$'s through
\begin{equation*}
  \epsilon_l = \frac{\epsilon}{\sqrt{\epsilon^2 + \left(\frac{\partial c}{\partial \theta}\Big|_{\theta=0}\right)^2}},
\end{equation*}
which, when applied to (\ref{descent}), yields descent steps of size bounded by $\epsilon$, and much smaller near the optimal $\theta$.

\bigskip
\bigskip

\begin{figure}[h]
\begin{center}
\vspace{-6.0cm}
\includegraphics[height=20cm]{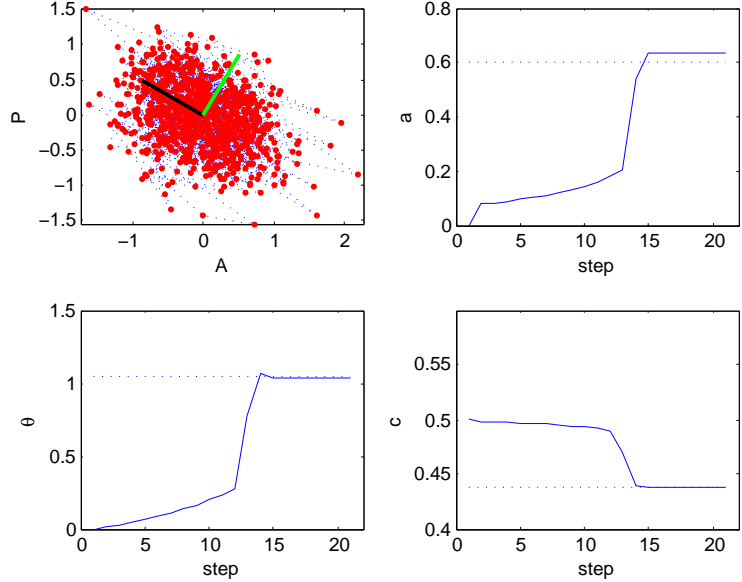}
\vspace{-7.0cm}
\end{center}
\caption{A two-dimensional example of the basic procedure. The first plot displays the data points, with dotted lines joining successive observations, and the directions for the first regular principal component --in black-- and the principal dynamical component --in green--, which in this case are orthogonal to each other. The other three plots show the evolution of the estimates for the parameters $a$ and $\theta$ for the dynamics and reduced manifold, and of and the cost function $c$  --normalized by $(N-1)$-- as functions of the step-number, with their exact values as dotted lines.}
\label{fig2Daut}
\end{figure}

To illustrate the procedure just described, we created data from the dynamical model
\begin{eqnarray*}
  x_{j+1} &=& a x_j +  r_x \eta^x_j, \\
  y_{j+1} &=&  r_y \eta^y_j ,
  \end{eqnarray*}
for $j=1, \ldots, N-1$ (the initial values $x_1$ and $y_1$ are picked at random), where $N=1000$, the $\eta^{x,y}_j$ are independent samples from a normal distribution, and we adopted the values $a=0.6$ for the dynamics\begin{footnote}{The value of $|a|$ needs to be smaller than one for the time series not to blow up.}\end{footnote}, and $r_x=0.3$ and $r_y=0.6$ for the amplitudes of the noise in $x$ and $y$.  Then we rotated the data through
\begin{eqnarray*}
A_j &=& x_j \cos(\theta_*) - y_j \sin(\theta_*), \\
P_j &=& x_j \sin(\theta_*) + y_j \cos(\theta_*) ,
\end{eqnarray*}
with $\theta_*=\frac{\pi}{3}$, and provided the $A_j$ and $P_j$ as data for the principal dynamical component routine. The results are displayed in Figure \ref{fig2Daut}. The first plot shows the ``observations'' in the plane $(A,P)$. These are treated as independent samples in a regular principal component analysis; we keep instead track of the sequential order of the observations, represented by the dotted lines in the plot. For this data, the first regular principal component, drawn in black, is in fact orthogonal to the principal dynamical component, drawn in green. The reason is that the total variability has a larger $y$-component, due to the bigger amplitude of the noise in $y$, while all the variability that is explainable dynamically is in $x$. This is an extreme example where regular principal components yield a leading mode that is absolutely irrelevant from a dynamical viewpoint. The other three plots in the figure display the evolution of the estimates for $a$ and $\theta_*$ and the cost function $c$, as functions of the step-number. The dotted lines, drawn for reference, have the exact values of $a$ and $\theta_*$ in the data, as well as the unexplainable part of the cost, $c_* = \frac{1}{N-1} \sum_{j=1}^{N-1} \left(r_x \eta^x_j\right)^2 +  \left(r_y \eta^y_j\right)^2 \approx r_x^2 + r_y^2 = 0.45$. Notice the fast convergence to the exact solution, that in this example took 14 steps.

\subsection{The multidimensional case}
\label{sec:linear}

For dimensions $n$ bigger than two, we write
\begin{equation*}
z = [Q_x Q_y] \left( \begin{array}{c}
		x \\
		y \end{array} \right) \, ,
\end{equation*}
and
\begin{equation*}
x_{j+1} = A x_j \, .
\end{equation*}
The minimization problem that defines $Q=[Q_x Q_y]$ and $A$ is
\begin{equation*}
\min_{Q,A}\, c = \sum_{j=1}^{N-1} \left\|z_{j+1} - Q \left(\begin{array}{c}
							A {Q_x}' z_j \\
							0 \end{array} \right) \right\|^2  =
		\sum_{j=1}^{N-1} \left\|\left(\begin{array}{c}
							x_{j+1} - A x_{j} \\
							y_{j+1} \end{array} \right)
							\right\|^2 \, .
\end{equation*}
Notice that $Q$ and $A$ are not univocally defined: any pair of orthogonal bases for the optimal subspaces represented by $x$ and $y$ will give rise to different $Q$'s and $A$'s representing the same dynamics. The algorithm proposed below walks nicely around this degeneracy, avoiding unnecessary re-parameterizations of the two subspaces.

The most straightforward methodology decouples the descent steps for $A$ and $Q$. For $A$, we have
\begin{equation}
\frac{\partial c}{\partial A} = -2 \sum_{j=1}^{N-1} \left(x_{j+1} - A x_j\right) x_j' \, .
\label{dcdA}
\end{equation}
%
%
Instead of descending the gradient we can, as in the two-dimensional case, directly solve $\frac{\partial c}{\partial A} = 0$, yielding
\begin{equation*}
A = X_1 {X_0}' \left(X_0 {X_0}'\right)^{-1} ,
\label{Aex}
\end{equation*}
where
\begin{equation*}
X_0 = \left[x_1, \ldots, x_{N-1}\right] \quad \hbox{and} \quad X_1 = \left[x_2, \ldots, x_N\right] \, .
\end{equation*}
For the descent steps in $Q$, we first note that any orthogonal matrix can be factorized
as a product of elementary rotations of the form
\medskip
\begin{equation*}
R_{kl}(\theta) =
\left( \begin{array}{ccccccc}
1 & \ldots & 0 & \ldots & 0 & \ldots & 0 \\
\vdots & \ddots & \vdots & \ddots & \vdots & \ddots & \vdots \\
0 & \ldots & \cos(\theta) & \ldots & \sin(\theta) & \ldots & 0 \\
\vdots & \ddots & \vdots & \ddots & \vdots & \ddots & \vdots \\
0 & \ldots & -\sin(\theta) & \ldots & \cos(\theta) & \ldots & 0 \\
\vdots & \ddots & \vdots & \ddots & \vdots & \ddots & \vdots \\
0 & \ldots & 0 & \ldots & 0 & \ldots & 1
		\end{array} \right) \, ,
\end{equation*}
which act in the plane of the two coordinates $k$ and $l$, rotating them an angle $\theta$.
Then we can, in each descent step, pick at random the two indices $k$ and $l$, and perform a rotation following the derivative of $c$ with respect to $\theta$ at $\theta=0$. From the observation above about degeneracy, however, we note that picking both $k$ and $l$ from either the dynamical coordinates $x$ or their orthogonal complement $y$ alone, serves no purpose other than re-parametrization. Then we always pick $k$ at random in $[1, \ldots, m]$, and adopt $l=m+h$, with $h$ picked at random in $[1, \ldots, n-m]$. In order to consider arbitrary directions in these two manifolds though, we first perform a random orthogonal transformation to each:
\begin{equation*}
  x \rightarrow Q_x^r x, \quad y \rightarrow Q_y^r y,
\end{equation*}
where $Q_{x,y}^r$ are random orthogonal matrices.

For each elementary rotation, we have
\begin{equation}
\frac{\partial c}{\partial \theta}\Big|_{\theta=0} = -2 \sum_{j=1}^{N-1} \left[ y_{j+1}^hA_k x_j+
y_j^h\sum_{p=1}^m A_p^k \left(x_{j+1}^p - A_p x_j \right) \right]
\label{ctheta}
\end{equation}
and
\begin{equation}
\frac{\partial^2 c}{\partial \theta^2}\Big|_{\theta=0}  = 2 \sum_{j=1}^{N-1} \left[ x_{j+1}^k A_kx_j -2y_{j+1}^hA_k^ky_j^h+
\sum_{p=1}^m \left[x_j^kA_p^k(x_{j+1}^p - A_p x_j)+(y_j^hA_p^k)^2\right] \right]\, .
\label{ctheta2}
\end{equation}
As before, $\theta$ can be computed so as to minimize the quadratic local approximation to $c$:
$$ \theta = \theta_q = -\frac{ \frac{\partial c}{\partial \theta}\Big|_{\theta=0} }{\frac{\partial^2 c}{\partial \theta^2}\Big|_{\theta=0} },$$
with the same caveats on big steps as in the one-dimensional case.

After performing the optimization, one can, if desired, resolve the degeneracy in the description of the dynamical manifold by choosing a natural basis for $x$, such as the one made out of the principal components of $A$. For non-normal $A$'s, there are two such bases: the eigenvectors of $A' A$ and those of $A A'$. Both are significant and sorted by
sensitivity to perturbations: the former gives the directions where initial perturbations yield the highest effect; the latter, the directions where these effects manifest themselves after a time-step.

\bigskip
\bigskip

\begin{figure}[h]
\begin{center}
\vspace{-6.0cm}
\includegraphics[height=20cm]{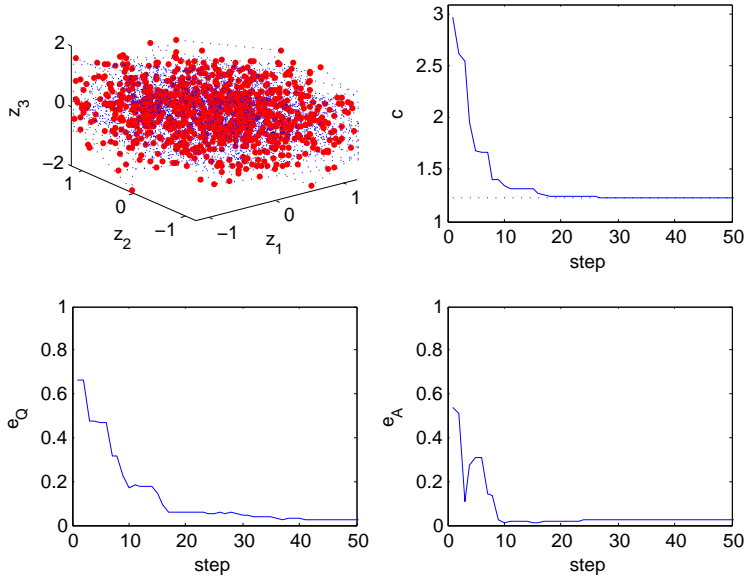}
\vspace{-7.0cm}
\end{center}
\caption{A multidimensional, autonomous example, with $n=5$ and $m=2$. The first plot displays the first three coordinates of the data points, with dotted lines joining successive observations, the second plot shows the evolution of the cost function, with its exact value as a dotted line, and the third and fourth plots display the evolution of the errors $e_Q$ and $e_A$.}
\label{figMDaut}
\end{figure}

To create a simple synthetic example for the multidimensional case, we chose $n=5$ and $m=2$, and created data from the dynamical model
\begin{eqnarray*}
  x_{j+1} &=& A x_j +  r_x \eta^x_j , \\
  y_{j+1} &=&  r_y \eta^y_j ,
  \end{eqnarray*}
for $j=1, \ldots, N-1$, where $N=1000$, the $\eta^{x,y}_j$'s are
independent samples from a normal distribution --two and three dimensional
vectors respectively-- and we adopted arbitrarily the values $$A=\begin{pmatrix}
0.4569   & 0.3237\\
   -1.0374  &  1.0378
\end{pmatrix}
$$ for the dynamics, and $r_x=0.3$ and $r_y=0.6$ for the amplitudes of the
noise in $x$ and $y$.  Then we rotated the data through an arbitrary orthogonal matrix,
\begin{equation*}
z = [Q_x Q_y] \left( \begin{array}{c}
                x \\
                y \end{array} \right) \, ,
\end{equation*}
with
\begin{equation*}
Q_x=\left( \begin{array}{ccccc}
   -0.7044  & -0.3823  & -0.3407  & -0.1985  & -0.4497\\
    0.5754  & -0.1555  & -0.1798  &  0.2477  & -0.7423
\end{array} \right)' \, ,
\end{equation*}
and generated the data displayed in the first panel of Figure \ref{figMDaut}. Running our algorithm on these data yields estimates $A^*$ and $Q_x^*$ for $A$ and $Q_x$ that, as remarked before, are not univocally defined. Indeed, the algorithm found
$$A^*=\begin{pmatrix}
    0.6505  &  0.2401\\
   -1.1591   & 0.8685
\end{pmatrix}
$$
and
\begin{equation*}
Q_x^*=\left( \begin{array}{ccccc}
   -0.8143  & -0.3258  & -0.2896 &  -0.2545 &  -0.2865\\
    0.4089  & -0.2108  & -0.2570 &   0.1873 &  -0.8290
\end{array} \right)' \, ,
\end{equation*}
quite different in appearance from their exact values above.

To verify that $Q_x$ and $Q_x^*$ span the same plane and that $A$ and $A^*$ represent the same transformation in the corresponding coordinates, we project the two columns of $Q_x^*$ onto the space spanned by those of $Q_x$, through the projection $P(Q_x^*) '= B Q_x'$, with
$B=\left(Q_x^*\right)'Q_x$, and define the relative errors
$$
e_Q=\frac{\|\left(Q_x^*\right)'-BQ_x'\|}{\|Q_x\|},\quad
e_A=\frac{\|A-B^{-1}A^*B\|}{\|A^*\|},
$$
which vanish only if the two pairs of matrices represent exactly the same reduced manifold and dynamics.

The results are displayed in Figure \ref{figMDaut}. The first plot shows the first three components of the data points  $z_j$. The second plot displays the evolution of the normalized cost function $c$ as a function of the step-number, with the dotted line displaying the exact value of the unexplainable part of the cost, $c_* = \frac{1}{N-1}
\sum_{j=1}^{N-1} \left(r_x \|\eta^x_j\|\right)^2 +  \left(r_y \|\eta^y_j\|\right)^2 \approx 2r_x^2 + 3r_y^2 = 1.26$. The third and fourth plots display the evolution of the errors $e_Q$ and $e_A$ defined above. Notice again the fast convergence of the algorithm to the exact solution.

\subsection{Non-zero means}

We have worked so far under the assumption that all means have been removed from the
problem: the plane $x$ goes through the origin, and the transformation given by the
matrix $A$ is linear, not affine. If the observations $z$ have a well-defined mean
(that is, if there is not a trend over time that makes the local mean of $z$ evolve),
these assumptions are fine: it is enough to remove from $z$ its mean --the ``climatology'' of atmosphere-ocean science-- ad initio, and add it back at the end. However, for the non-autonomous scenario to be described below, it will be necessary to consider nontrivial means. In order to have our methodology prepared for this more general case, we consider the means in our present autonomous situation too, even though they have no practical consequence. Then we write
\begin{equation*}
z - \bar{z} = Q \left( \begin{array}{c}
		x \\
		y \end{array} \right) \, ,
\end{equation*}
and
\begin{equation*}
x_{j+1} = A x_j + b \, .
\end{equation*}
It is convenient to partition $\bar{z}$ into its $x$ and $y$ components,
\begin{equation*}
\bar{x} = Q_x' \bar{z} \, , \quad  \bar{y} = Q_y' \bar{z} .
\end{equation*}
The addition of the mean $\bar{x}$, however, is unnecessary, for its effects can be absorbed into the drift $b$. We have the gradients
\begin{equation}
\frac{\partial c}{\partial b} = -2 \sum_{j=1}^{N-1} \left(x_{j+1} - (A x_j + b)\right)
\label{dcdb}
\end{equation}
and
\begin{equation}
\frac{\partial c}{\partial \bar{y}} = -2 \sum_{j=1}^{N-1} \left(y_{j+1}-\bar{y}\right) \, ,
\label{dcdy}
\end{equation}
that can be used either for descent or for the direct calculation of the optimal $b$ and $\bar{y}$.

\section{Non-autonomous problems}
\label{NonAutonomous}

We have considered up to now only autonomous problems, where the manifold $x$ and the corresponding dynamical model are assumed to be time-independent. Yet there are many examples of practical importance where this assumption does not hold. Consider, for instance, climate-related data, such as monthly averages of sea-surface temperatures at various locations, recorded over many years. One should expect much of the dynamics to depend on seasonal changes in insolation. We should, accordingly, have a time-dependent dynamical model, with a period of one year. Similarly, in long series of economic or financial data, we should expect a change in the dynamics as populations or affluence levels change, new markets arise, new tools are developed. The corresponding dynamical model should not longer be constant, nor periodic as in the seasonal case, but rather evolve slowly, with scale separation between the time-scale of the dynamics and that of the evolution of the model itself (without the hypothesis of scale separation, little can be inferred statistically from the data, since the dynamical model can be adjusted instantly to account for each individual observation).

To incorporate this into our framework, it is enough to add a qualifying sub-index ``$t$'' (or more precisely ``$j$'', since our time-series are discrete) to the various functions involved: $Q$, $A$, $b$ and $\bar{y}$ , plus the requirements of periodicity or scale separation. For instance, $Q_t$ should satisfy either $Q_{t+T} = Q_t$ in the periodic case, or $\left\|Q_{t+1}-Q_t \right\| \ll 1$ for slowly varying trends. In this section, we discuss how to modify the methodology of Section \ref{sec:linear} so as to make it applicable to the non-autonomous linear case.

The idea is simple: in the notation of the previous section, we are seeking a time-dependent orthogonal transformation $Q_t$ and mean $\bar{y}_t$, and a time-dependent dynamical model parameterized by $A_t$ and $b_t$. To this end, in each descent step, we pick at random a time $t_0$ and propose, in order to update $Q$, a time-dependent rotation angle $\theta(t)$ in the $k$-$l$ plane, centered at $t=t_0$. Similarly,  we propose time-dependent variations for  $\bar{y}$, $A$ and $b$:
\begin{equation*}
\theta = \alpha F(t) \, , \quad \bar{y} = \bar{y} + v \, F(t)\, ,\quad A = A + B \, F(t) \, ,\quad b = b + d \, F(t) \, ,
\end{equation*}
where $F(t)$ is a given scalar function, centered at $t_0$, and satisfying the corresponding restrictions: periodicity, slow variation, etc., and the parameters $\alpha$, a scalar, $v$ and $d$, vectors, and $B$, a matrix, are computed by descent of the cost function as before. Then equations (\ref{ctheta}) and (\ref{ctheta2}) generalize into
\begin{equation*}
\frac{\partial c}{\partial \alpha}\Big|_{\alpha=0} = -2 \sum_{j=1}^{N-1} \left[ w^{j+1}y_{j+1}^h\left(A_k x_j+b_k\right)+
w^jy_j^h\sum_{p=1}^m A_p^k \left(x_{j+1}^p - A_p x_j-b_p \right) \right]
\label{calpha}
\end{equation*}
and
\begin{align*}\label{calpha2}
\nonumber\frac{\partial^2 c}{\partial \alpha^2}\Big|_{\alpha=0} &= 2 \sum_{j=1}^{N-1} \bigg[ (w^{j+1})^2x_{j+1}^k \left(A_kx_j+b_k\right) -2w^{j+1}w^{j}y_{j+1}^hA_k^ky_j^h+\\
 &\qquad \qquad \sum_{p=1}^m \left[(w^{j})^2x_j^kA_p^k(x_{j+1}^p - A_p x_j-b_p)+(w^{j}y_j^hA_p^k)^2\right] \bigg]\, ,
\end{align*}
and equations (\ref{dcdA}), (\ref{dcdb}) and (\ref{dcdy}) into
\begin{equation*}
\frac{\partial c}{\partial B} = -2 \sum_{j=1}^{N-1} w^j \left(x_{j+1} - (A x_j+b)\right) x_j' \, ,
\end{equation*}
\begin{equation*}
\frac{\partial c}{\partial d} = -2 \sum_{j=1}^{N-1} w^j\left(x_{j+1} - (A x_j + b)\right)
\label{dcdd}
\end{equation*}
and
\begin{equation*}
\frac{\partial c}{\partial v} = -2 \sum_{j=1}^{N-1} w^j \left(y_{j+1}-\bar{y}\right) \, ,
\label{dcdv}
\end{equation*}
where the weights $w^j$ are given by
\begin{equation*}
w^j = F(t_j) \, .
\end{equation*}
As before, equating the derivatives with respect to $B$, $d$ and $v$ to zero provides simple closed forms for $B$, $d$ and $v$, while $\alpha$ can be found through the minimization of a local quadratic approximation to $c$.

\subsection{Exogenous variables}

The time $t$ of the non-autonomous scenario discussed above is just one example
of an exogenous variable: one whose state is known independently at all times, and that may affect the dynamics of the $z$'s. Other examples are state variables of a bigger system of which the $z$'s are only a small part; and external controls.

One can collectively denote these exogenous variables $s$, and apply a straightforward
generalization of the procedure above, where $F$ is now a function of $s$ rather than the single variable $t$.

\subsection{Trial functions}
\label{trials}

We have not yet considered the issue of how to pick the functions $F(s)$ and
corresponding weights $w^j = F(s^j)$ (here we use $s$ to denote either time or other exogenous variables). In this section, we describe a few choices that we have found practical. First of all, for the autonomous case, we have the trivial
$$ F = 1 \, . $$
This should still be used in the more general case, to capture the $s$-independent components of $Q$ and $A$, but must be alternated with other functions $F(s)$ with non-trivial $s$-dependence.

\subsubsection{One-dimensional functions}

\begin{figure}[h]
\begin{center}
\vspace{-4.9cm}
\includegraphics[height=15cm]{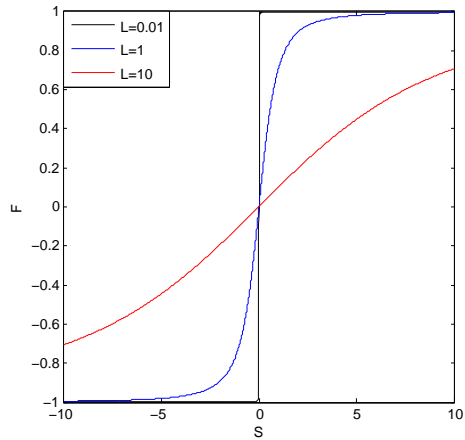}
\vspace{-5.2cm}
\end{center}
\caption{Plot of the function $F(S)=\frac{S}{\sqrt{S^2+L^2}}$ over the interval $[-10,10]$ for different values of $L$.}
\label{figTT1}
\end{figure}

When $s$ represents time, the domain of $F(s)$ must be either the real line --for the trend-- or a parametrization of the unit circle --for periodic factors such as the seasons.
A sensible choice for the trend is
\begin{equation*}
 F(s) = \frac{S}{\sqrt{S^2+L^2}} - \bar{F} , \quad S=s-s_0,
 \label{trend}
\end{equation*}
displayed on Figure \ref{figTT1}, depending on the choice of a center $s_0$, picked at random at each step, and a mollification parameter, the length-scale $L$.  As $L \rightarrow 0$, $F(s)$ becomes piecewise constant, with a discontinuity at $s=s_0$. For larger values of $L$, the transition between the two constant states is smoothed over an interval of order $L$. The subtraction of the mean $\bar{F}$ over the observations is intended to decouple the
effect of these steps from the ones using $F=1$, a function concerned only with the mean.
At the initial stages of the algorithm, $L$ should be large, providing a global perspective;
then it should decrease gradually, to tune the finer, more local details.

We can be more specific: calling $L_0$ the largest length scale in $s$, we need $L_0/L$ steps to cover it with transitions of length $L$. Then the amount $dt$ of algorithmic time
spent using a length $L$ should satisfy
$$ \frac{dL}{dt} \propto L ,$$
leading to the expression
$$ L=L_0 \left(\frac{L_f}{L_0}\right)^{\frac{k}{k_{tot}}} \, , $$
where $k$ is the step number, $k_{tot}$ the total number of steps, and $L_f$
the smallest length scale to be used, not to over-resolve the dynamics.

\begin{figure}[h]
\begin{center}
\vspace{-4.9cm}
\includegraphics[height=15cm]{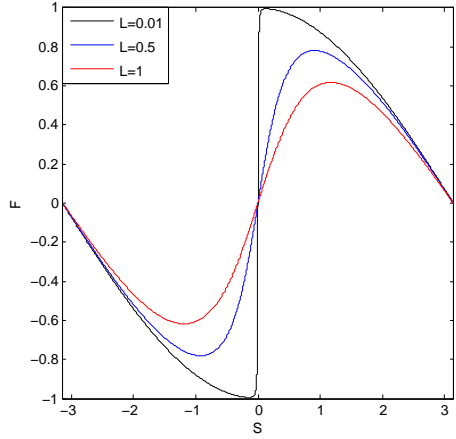}
\vspace{-5.2cm}
\end{center}
\caption{Plot of the function $F(S)=\frac{\sin(S)}{\sqrt{4 \sin^2(S/2) + L^2}}$ over the interval $[-\pi,\pi]$ for different values of $L$.}
\label{figTP1}
\end{figure}

In the periodic case, we can make an entirely analogous proposal:
\begin{equation}\label{Wfs}
  F(s) = \frac{\sin(S)}{\sqrt{4 \sin^2(S/2) + L^2}}, \quad S = \frac{2\pi \left(s-s_0\right)}{T}\, ,
\end{equation}
where $T$ is the period; see Figure \ref{figTP1}.

Sometimes $s$ can adopt only a discrete set of values: the months of the year,
an on-off control, etc. In that case, it may be useful to consider signature functions $F$
that are one on each of these values at a time, and zero on the others:
\begin{equation}
  F_i(j) = \delta_{mod(j,T),i} ,
  \label{delta}
\end{equation}
where $T$ is the integer period.

\begin{figure}[h]
\begin{center}
\vspace{-4.9cm}
\includegraphics[height=15cm]{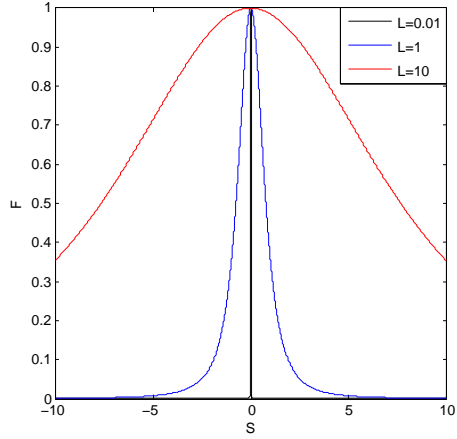}
\vspace{-5.2cm}
\end{center}
\caption{Plot of the function $F(S)=\frac{L^3}{(S^2+L^2)^{3/2}}$ over the interval $[-10,10]$ for different values of $L$.}
\label{figTT2}
\end{figure}

\begin{figure}[h]
\begin{center}
\vspace{-4.9cm}
\includegraphics[height=15cm]{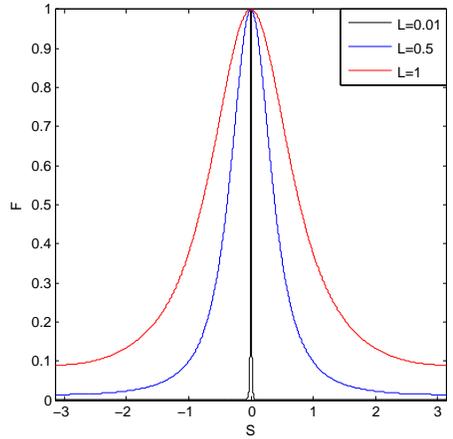}
\vspace{-5.2cm}
\end{center}
\caption{Plot of the function $F(S)=\frac{L^3}{(4 \sin^2(S/2) + L^2)^{3/2}}$ over the interval $[-\pi,\pi]$ for different values of $L$.}
\label{figTP2}
\end{figure}

An alternative to the $F(s)$'s above, which have local derivatives but global
effects, are the more localized bumps given by
$$ F(s) = \frac{L^3}{(S^2+L^2)^{3/2}} - \bar{F} , $$
and
$$ F(s) = \frac{L^3}{(4 \sin^2(S/2) + L^2)^{3/2}} - \bar{F}  ,$$
displayed in Figures \ref{figTT2} and \ref{figTP2}. 
We can also alternate between the two, or among more proposals satisfying
different needs.


Still another natural alternative in the periodic case is to use Fourier components
$$ F_k^c = \cos(k S) , \quad F_k^s = \sin(k S) , \quad S = \frac{2\pi s}{T}. $$
There is no need for a center $s_0$ here, since the use of both sines and cosines renders the $F$ spatially homogeneous. Each step one must use either $F^c$ or $F^s$ with probability $1/2$ each (discounting the steps with $F=1$), and an integer value for the wave number $k$. The latter should be sampled from a distribution that decays rapidly with $k$, so as to result into smooth composite functions.

An advantage of the use of Fourier modes, particularly when only a finite number $K$ of modes is allowed, is that one can store the accumulated amplitude added to each mode at the various steps, and thus end up with explicit expressions for the non-autonomous dynamical matrix $A(s)$ and shift $b(s)$ as finite Fourier series. To obtain a similar bonus for the non-periodic case (i.e., for representations of the trend), one would need to replace the functions above by others that do not involve a variable length-scale $L$ and random point $s_0$. A simple choice is that of monomials
$$ F_k(s) = s^k ,$$
up to a power $K$. Then the dynamics is represented by a matrix $A$ and a vector $b$ that depend explicitly on $s$ through polynomials of degree $K$.

When the problem has more that one kind of variable --some periodic and some trendy, for instance--, we can alternate the various types of function $F(s)$ among steps.

\subsubsection{Multidimensional choices}

When $s$ lives in a multidimensional space, we can still use the one-dimensional proposals involving $s_0$ and $L$ above, but picking the direction of space in which they apply each step at random. Yet this is not a very effective procedure when the dimensionality of $s$ is large. An alternative is to use radial functions centered at $s_0$, such as
%
%
$$ F(s) = e^{-r^2}, \quad r = \frac{\|s-s_0\|}{L} .$$

\subsection{A non-autonomous example}

For clarity, we illustrate the non-autonomous procedure through a simple example where $n=2, m=1$. We created data from the dynamical model
\begin{eqnarray*}
  x_{j+1} &=& a_jx_j +b_j+  r_x \eta^x_j ,\\
  y_{j+1} &=&  \bar{y}_{j+1} + r_y \eta^y_j ,
  \end{eqnarray*}
for $j=1, \ldots, N-1$, where $N=1000$, the $\eta^{x,y}_j$'s are independent samples from a normal distribution, $r_x=0.3$ and $r_y=0.6$, and we adopted the values
$a_j=\frac{6}{5}\cos^2\left(\frac{2\pi t_j}{T}\right)$ for the dynamics,
$b_j=\frac{1}{2}\sin\left(\frac{2\pi t_j}{T}\right)$ for the drift, and
$\bar{y}_j=\frac{2}{5}\cos\left(\frac{2\pi t_j}{T}\right)$ for the
non-zero mean of $y$, where $t_j=j$ and $T=12$, mimicking the twelve months of the year that we will find again in our application to the sea-surface temperature field in Section \ref{Realex}. Then we introduce, as before, ``Atlantic'' and ``Pacific'' temperatures
\begin{eqnarray*}
A_j &=& x_j \cos(\theta_j) - y_j \sin(\theta_j) ,\\
P_j &=& x_j \sin(\theta_j) + y_j \cos(\theta_j) ,
\end{eqnarray*}
with $\theta_j=\frac{\pi}{6}\sin\left(\frac{2\pi t_j}{T}\right)$, and provide the $A_j$ and $P_j$ as data for the principal dynamical component routine.

\begin{figure}[h]
\begin{center}
\vspace{-6.0cm}
\includegraphics[height=20cm]{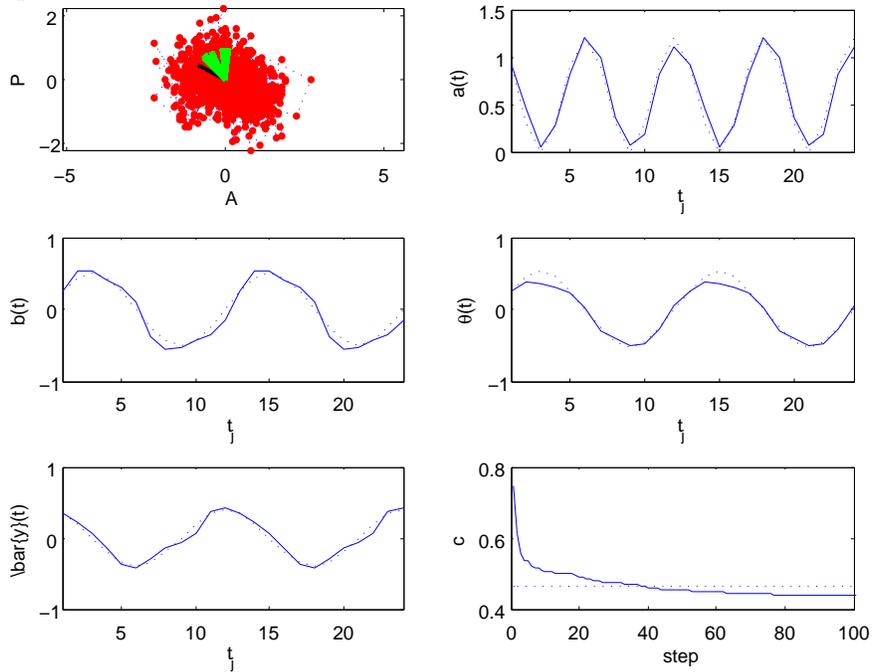}
\vspace{-6.5cm}
\end{center}
\caption{A low dimensional ($n=2, m=1$), non-autonomous problem. The first plot displays the data points, with dotted lines joining successive observations, the first principal component in black, and the twelve monthly first principal dynamic components in green. The other plots show  the final estimates for $a(t)$, $b(t)$, $\bar{y}(t)$ and $\theta(t)$, for two periods of twelve snap-shots each, and the evolution of the cost function, with the exact answers in dotted lines.}
\label{fig2DNonaut}
\end{figure}

For this example, we have adopted the trial function $F$ from (\ref{Wfs}). The results are displayed in Figure \ref{fig2DNonaut}. The first plot shows the ``observations'' in the plane $(A,P)$, with the first regular principal component drawn in black, and the 12 first principal dynamical components, one for each month, drawn in green.
The other plots in the figure display the evolution of the normalized cost function $c$, and the estimated results for $a(t)$, $b(t)$, $\bar{y}(t)$ and $\theta(t)$ at convergence (we only show the first two periods). The dotted lines, drawn for reference, have the exact values of $a(t)$, $b(t)$, $\bar{y}(t)$ and $\theta(t)$ in the data, as well as the unexplainable part of the cost, $c_* = \frac{1}{N-1} \sum_{j=1}^{N-1} \left(r_x \eta^x_j\right)^2 +  \left(r_y \eta^y_j\right)^2 \approx r_x^2 + r_y^2 = 0.45$. Again, the algorithm detects essentially the exact solution to the problem; the number of required steps, about 60, is bigger than before, because various different trial functions $F(t)$ are involved, requiring at least one step for each.

\section{Higher order processes}
\label{NonM}

We have considered so far dynamical models without memory, where the current state
of the system determines its future evolution through the matrix $A$.
Yet many real processes are not well-described by such models. For instance, if the
observations consist only of positions $x_j$ in a system with non-negligible inertia, one
would expect a better prediction by using, in lieu of the unavailable velocity field,
a second order model, $x_{j+1} = D(x_j, x_{j-1})$.
Studying systems like this involves no significant change in our procedure: either we extend the phase-space from the line of $x_j$'s to the plane of pairs $(x_j, x_{j-1})$ or, equivalently, consider matrices $A$ that are rectangular, with twice as many columns as rows. Entirely similar considerations apply to higher order processes with longer memory.

We describe here the non-Markovian, non-autonomous case of order $r$, since the autonomous scenario is just a special case of the non-autonomous one, and the case with more general exogenous variables $s$ is entirely similar. Our reduced dynamical model now adopts the form
\begin{equation*}
x_{j+1}=D=b+\sum_{i=1}^{r}A_i x_{j-i+1},
\end{equation*}
where the drift $b$ and the matrices $A_i, i=1,\ldots, r$, as well as the orthogonal matrix $Q_x$ defining the $x$'s, may in general be time-dependent. Each algorithmic step, we update these matrices through
$$
A_i=A_i+B_i\, F(t)\, , \quad b=b+d\, F(t)\, ,\quad\bar{y} = \bar{y} + v \, F(t)\, ,\quad \theta=\alpha \, F(t)\, ,
$$
where $F(t)$ is a given trial function as described above.

The cost function adopts the form
\begin{equation*}\label{costNM}
    c=\sum_{j=r}^{N-1} \left\|y_{j+1}\right\|^2 + \left\|x_{j+1}-D\right\|^2,
\end{equation*}
since the first $x_1,\ldots,x_r$ are not specified by the dynamics.
Then
\begin{equation*}
\frac{\partial c}{\partial B_h} = -2 \sum_{j=r}^{N-1} w^j\left(x_{j+1} - D\right) x_{j-h+1}' \, ,
\label{dcdBkNM}
\end{equation*}
for every $h=1,\ldots,r$,
\begin{equation*}
\frac{\partial c}{\partial d} = -2 \sum_{j=r}^{N-1} w^j\left(x_{j+1} - D\right) \, ,
\label{dcddNM}
\end{equation*}
and
\begin{equation*}
\frac{\partial c}{\partial v} = -2 \sum_{j=r}^{N-1} w^j \left(y_{j+1}-\bar{y}\right) \, ,
\label{dcdvNM}
\end{equation*}
where $w^j=F(t_j)$.

It is possible to get explicit expressions for $B_1,\ldots,B_r$ by equating all $\frac{\partial c}{\partial B_h}$  to zero.  We introduce the $m\times m$ matrices
$$ X_0^{h,k}=\sum_{j=r}^{N-1}\left(w^j\right)^2x_{j-h+1}x_{j-k+1}',\quad X_1^{h}=\sum_{j=r}^{N-1}w^j\left(x_{j+1}-b\right)x_{j-h+1}',
$$
where $X_0^{h,k}=\left(X_0^{k,h}\right)'$.
In terms of these, we get a block system of linear equations:
\begin{align*}
 B_1X_0^{1,1}+ B_2X_0^{2,1}+\cdots+B_rX_0^{r,1}=&X_1^{1},\\
 B_1X_0^{1,2}+ B_2X_0^{2,2}+\cdots+B_rX_0^{r,2}=&X_1^{2},\\
\vdots\qquad\qquad\vdots\qquad\qquad\qquad\vdots\qquad\qquad &\;\;\vdots\\
B_1X_0^{1,r}+ B_2X_0^{2,r}+\cdots+B_rX_0^{r,r}=&X_1^{r},
\end{align*}
or
\begin{equation*}
 \begin{bmatrix}
B_1, & B_2, & \cdots & B_r
\end{bmatrix}
\begin{bmatrix}
X_0^{1,1} & X_0^{1,2} & \cdots & X_0^{1,r}\\
X_0^{2,1} & X_0^{2,2} & \cdots & X_0^{2,r}\\
\vdots&\vdots&\ddots&\vdots\\
X_0^{r,1} & X_0^{r,2} & \cdots & X_0^{r,r}\\
\end{bmatrix}= \begin{bmatrix}
X_1^1, & X_1^2, & \cdots & X_1^r
\end{bmatrix},
\end{equation*}
which determines the matrices $B_1,\ldots,B_r$.

For the angle $\alpha$ we proceed as in the previous sections, though a quadratic approximation to $c$, using
\begin{equation*}
\frac{\partial c}{\partial \alpha}\Big|_{\alpha=0}=-2 \sum_{j=r}^{N-1}\bigg[w^{j+1}y_{j+1}^h\left(D\right)_k+\sum_{p=1}^m (D_{\alpha})_p\left(x_{j+1}^p -\left(D\right)_p \right) \bigg],
\end{equation*}
and
\begin{align*}
\frac{\partial^2 c}{\partial \alpha^2}\Big|_{\alpha=0}=& -2 \sum_{j=r}^{N-1} \bigg[ -\left(w^{j+1}\right)^2x_{j+1}^k \left(D\right)_k +2w^{j+1}y_{j+1}^h(D_{\alpha})_k\\
&\qquad\qquad\qquad+\sum_{p=1}^m \left[(D_{\alpha\alpha})_p\left(x_{j+1}^p-\left(D\right)_p\right)-\left[(D_{\alpha})_p\right]^2\right] \bigg],
\end{align*}
where
\begin{equation*}
    (D_{\alpha})=\sum_{i=1}^{r}(A_i)^kw^{j-i+1}y_{j-i+1}^h, \quad
    (D_{\alpha\alpha})=-\sum_{i=1}^{r}(A_i)^k\left(w^{j-i+1}\right)^2x_{j-i+1}^k.
\end{equation*}

\bigskip
\bigskip

\begin{figure}[h]
\begin{center}
\vspace{-6.0cm}
\includegraphics[height=20cm]{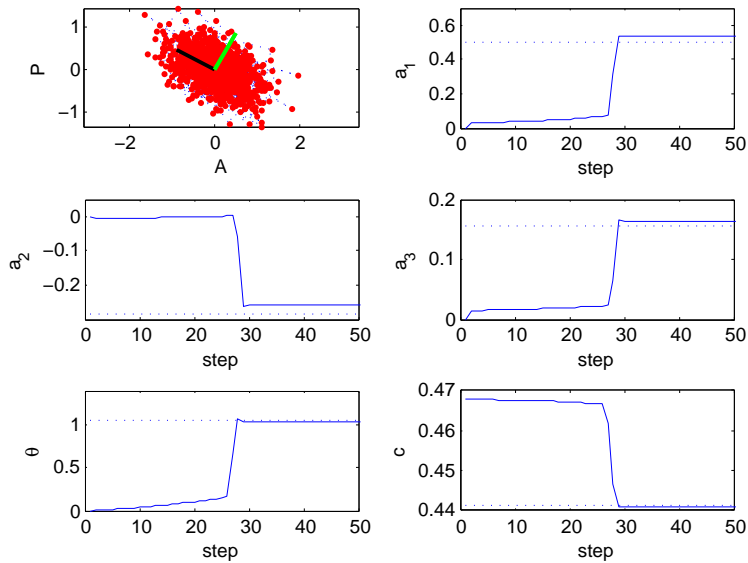}
\vspace{-7.0cm}
\end{center}
\caption{A multi-step process of order 3. The first plot displays the data points, with dotted lines joining successive observations, their regular first principal component in black and their first principal dynamical component in green. The other plots show the evolution of the estimates for $a_1,a_2,a_3$ and $\theta$, as well as of the cost function, with their exact values and the exact unexplainable part of the cost, $c_* = \frac{1}{N-1} \sum_{j=1}^{N-1} \left(r_x \eta^x_j\right)^2 +  \left(r_y \eta^y_j\right)^2 \approx r_x^2 + r_y^2 = 0.45$, displayed in dotted lines.}
\label{fig2DNonMar}
\end{figure}

Again we choose, for the sake of clarity, to illustrate the procedure in its simplest possible setting, which is autonomous, with $n=2$, $m=1$, and $r=3$, the order of the Non-Markovian process. We created data from the dynamical model
\begin{eqnarray*}
  x_{j+1} &=& a_1x_j+a_2x_{j-1}+a_3x_{j-2}+  r_x \eta^x_j ,\\
  y_{j+1} &=&  r_y \eta^y_j ,
  \end{eqnarray*}
for $j=3, \ldots, N-1$, where $N=1000$, the $\eta^{x,y}_j$'s are independent samples from a normal distribution, and we adopted the values $a_1=0.4979, a_2=-0.2846, a_3=0.1569$ for the dynamics and $r_x=0.3$ and $r_y=0.6$ for the amplitudes of the noise in $x$ and $y$. Then, as before, we define
\begin{eqnarray*}
A_j &=& x_j \cos(\theta_*) - y_j \sin(\theta_*) ,\\
P_j &=& x_j \sin(\theta_*) + y_j \cos(\theta_*) ,
\end{eqnarray*}
with $\theta_*=\frac{\pi}{3}$, and provide the $A_j$ and $P_j$ as data for the principal dynamical component routine. The results are displayed in Figure \ref{fig2DNonMar}. Again the procedure converges to the exact answer, this time for all elements of the multi-step dynamics. As in the first example,  the first regular principal component is orthogonal to the principal dynamical component, thus capturing none of the system's dynamics.

\section{A real application: the global sea-surface temperature field}
\label{Realex}

To see the workings of the new procedure on real data, we have chosen a topic of present concern: the estimation of climatic variations and trends. For this, we use a database of monthly averaged extended reconstructed global sea surface temperatures based on COADS data (see \cite{IRI}) from January 1854 to October 2009, and ask whether we can extract from these a reduced low dimensional dynamical model. A few before-hand considerations are in order:
\begin{itemize}
  \item Climate dynamics, a real pressing issue, is treated just as an illustration in this methodological paper. A far more in-depth treatment of how much principal dynamics components can help increase climate predictability and elucidate its causal relations will be pursued elsewhere.

  \item The ocean is not an isolated player in climate dynamics: it interacts with the atmosphere and the continents, and is also affected by external conditions, such as interannual variations in solar radiation and human-related release of $CO_2$ into the atmosphere. The latter are examples of slowly varying external trends that fit naturally into our non-autonomous setting --the seasonal variations giving its periodic component. As for the land and atmosphere, their dynamics is typically faster than that of the oceans, and can be conceptually divided into two components: a part that is slaved to the state of the ocean's surface temperature --and hence can in principle be included in its dynamical model--, and one that can be treated as external noise. Including explicitly land and atmospheric observations involves at least two further challenges, that will be pursued elsewhere: handling data with disparate units --such as atmospheric pressure, ice extent and ocean temperature--, and allowing for multiple time-scale dynamical models.

  \item Even within the ocean, the surface temperature does not evolve alone: it is carried by currents, and it interacts through mixing with lower layers of the ocean. As mentioned in Section \ref{NonM}, one way to account for unobserved variables is to make the model non-Markovian: discrete time derivatives of the sea-surface temperature provide indirect evidence on the state of those hidden variables.

\end{itemize}

\begin{figure}[h]
\begin{center}
\vspace{-5.0cm}
\includegraphics[height=20cm]{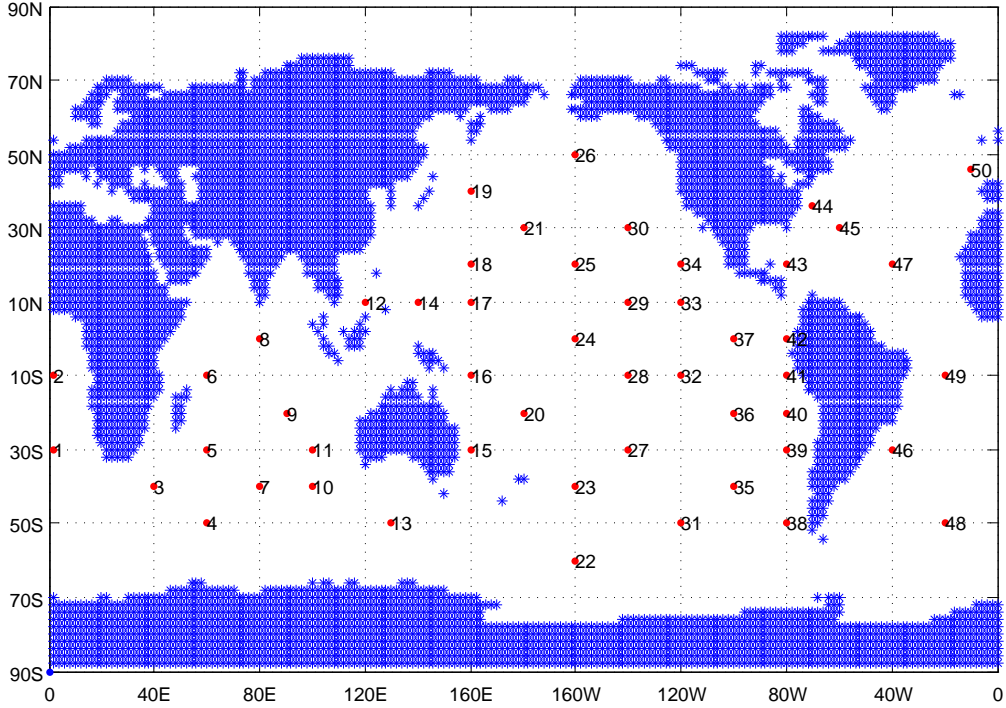}
\vspace{-6.0cm}
\end{center}
\caption{The 50 points on the ocean used for the procedure.}
\label{fig50points}
\end{figure}

We have adopted as our dataset the sea-surface temperature monthly means between January 1854 to October 2009 of the 50 points displayed on the map in Figure \ref{fig50points}, covering much of the world oceans in a roughly homogeneous manner.

In order to apply our methodology to the data, we need to select a class of trial functions from Subsection \ref{trials}, the dimension $m$ for the reduced manifold $x$, and the order $r$ of the non-Markovian process. The trial functions for the periodic component that we have used for these runs are the monthly discrete $\delta$-functions  from (\ref{delta}), with $T=12$. This takes to a new depth the idea behind the use of a ``monthly climatology'' in climate studies: not only the climatological mean is computed independently for each month, but also the dynamical model and manifold may change significantly from month to month. In the runs reported here, we have not modeled any inter-annual trend.

The following physical considerations suggest picking $r=3$ for the order of the Markov process. A simplified conceptual model for the upper mixed layer of the ocean is that of a rotating shallow layer of water, forced by the atmosphere from above and the deep ocean from below. In such model, the active dynamical variables are the two horizontal components of the velocity and the layer's thickness. The surface temperature can be thought of as an emergent of the evolution of these three variables and the external forcing. Conservation of mass and horizontal momentum, the core dynamics of the layer, are three differential equations, each involving one time derivative. Hence reducing the system to a single variable --the temperature, the only one available in the data-- yields a third order differential equation: two time derivatives relate to the evolution of gravity waves, the third to the potential vorticity. In our discrete setting, this corresponds to a Markov process of order three.

\begin{figure}[h]
\begin{center}
\vspace{-7.0cm}
\includegraphics[height=20cm]{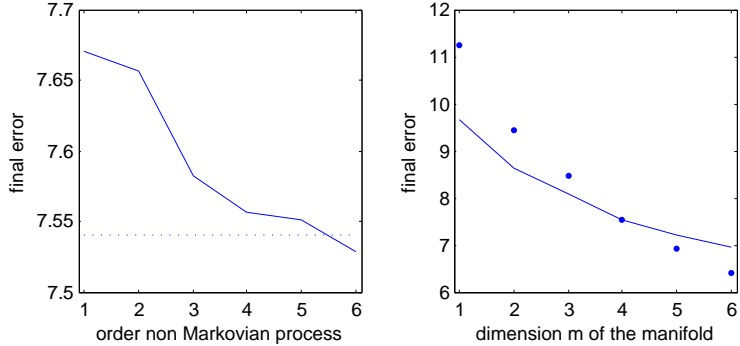}
\vspace{-8.0cm}
\end{center}
\caption{Predictive uncertainty as a function of the dimension $m$ of the reduced dynamical manifold and the order $r$ of the process. The figure on the left shows the final error for a fixed $m=4$ and $r=1,\ldots,6$. The dotted line is the sum $\frac{1}{N}\sum_{i=5}^NS_i^2$ of the squared singular values of the data, with the monthly climatology removed. The figure on the left shows the final error for a fixed $r=3$ and $m=1,\ldots,6$, with the isolated points displaying the corresponding uncertainty in the standard principal component procedure,  $\frac{1}{N}\sum_{i=m+1}^NS_i^2$.}
\label{figorderdimen}
\end{figure}

Figure \ref{figorderdimen} illustrates a line of reasoning for choosing the values of $m$ and $r$. The figure on the left shows, for a fixed $m=4$, the evolution of the final error when we move the order of the non Markovian process from $r=1$ to $r=6$\begin{footnote}{This final error is calculated for a number of steps such that the difference between the final error and the error 1000 steps before is less than 0.01. Therefore the number of steps used for each value of $r$ may be different.}\end{footnote}.  The dotted line shows the error $\frac{1}{N}\sum_{i=5}^NS_i^2$, where the $S$'s are the singular values of the real dataset, with the monthly climatology subtracted. We find for $r=3$ the steepest drop of the final error, consistent with our reasoning above. Therefore we pick $r=3$ for the order of our non Markovian process. In the figure on the right, we observe, for this fixed $r=3$, the evolution of the final error when $m=1,\ldots,6$. The isolated points correspond to the sum of squared singular values, $\frac{1}{N}\sum_{i=m+1}^NS_i^2$. We observe that for $m=4$ this error matches almost exactly the one from the dynamical components. This can be interpreted in the following way: for smaller values of $m$, accounting for the dynamics allows us to reduce the information loss even beyond the theoretical maximum --for autonomous settings-- provided by the singular value decomposition. Beyond $m=4$, on the other hand, the biggest share in the further reduction of information loss is probably due to the increased bare dimensionality of the model, more than to a further refinement of the dynamics. Hence we pick $m=4$ for the dimension of our reduced dynamical manifold.

\bigskip

\begin{figure}[h]
\begin{center}
\vspace{-6.0cm}
\includegraphics[height=20cm]{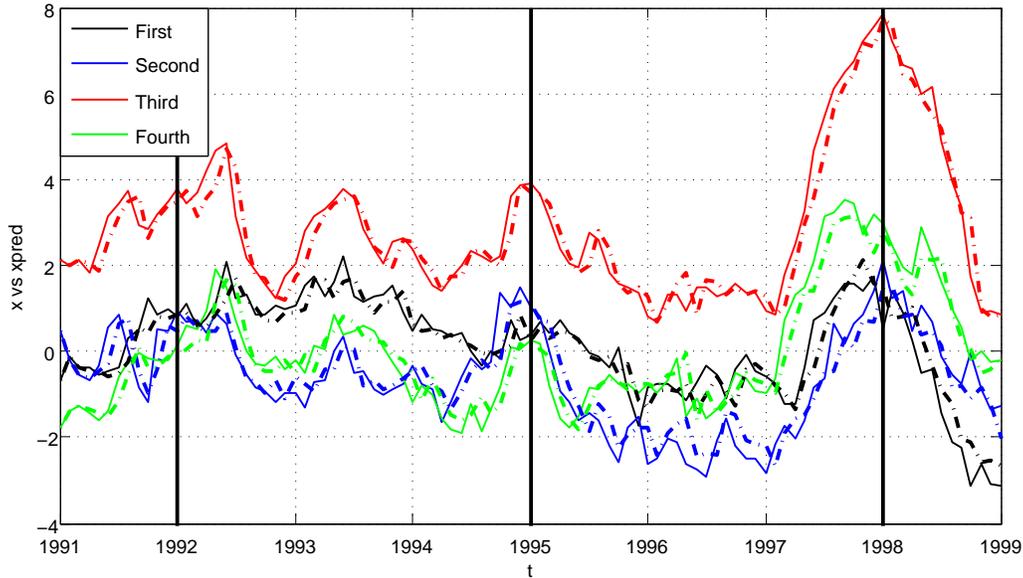}
\vspace{-7.0cm}
\end{center}
\caption{Real and predicted dynamical components (the manifold $x$). The vertical lines mark El Ni\~no years.}
\label{xvsxpred}
\end{figure}

Next we display various results for the chosen parameters, $m=4$ and $r=3$. Although  the dataset used is defined over a wide range of time, we display our results on the time window from January 1991 to January 1999, which includes three El Ni\~{n}o years, represented by vertical lines; one of them, in 1998, the strongest ever recorded.
Figure \ref{xvsxpred} shows the evolution of the four components of the manifold $x$ in solid lines and, in dotted lines, the same components predicted from the prior three months.

\begin{figure}[h]
\begin{center}
\vspace{-6.0cm}
\includegraphics[height=20cm]{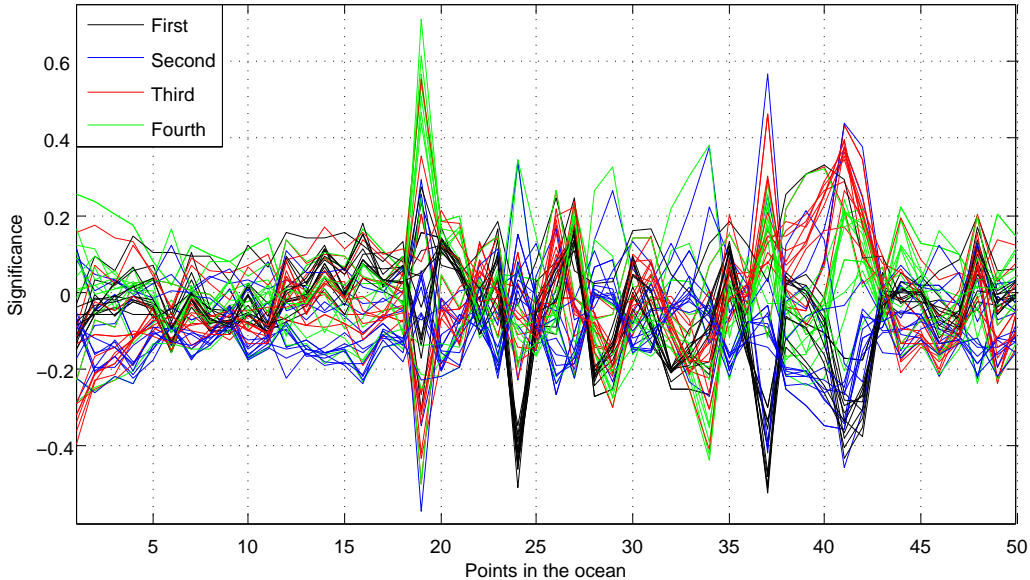}
\vspace{-7.0cm}
\end{center}
\caption{Dependence of the four components of the reduced manifold $x$ on the individual locations on the ocean for the twelve months of the year. A natural coordinate system  in $x$ is the one provided by the principal components of the first dynamical matrix $A_1$. Notice that four to six points on the ocean dominate the dynamics.}
\label{SigPoints}
\end{figure}

One question one may ask is whether the reduced dynamical manifold $x$ is dominated by a small set of locations on the ocean. This would be manifested in having the columns of $Q_x$ dominated by a few significant rows. Yet the columns of $Q_x$ do not have a meaning per se: $x$ is a four dimensional manifold, but each component $x^i$ lacks individual meaning. To fix a reference frame in the manifold $x$, we resort to the matrix $A_1$: its four left principal components $U$ in $A_1 = U  S  V'$ provide a natural set of coordinates in $x$-space.
Figure \ref{SigPoints} displays the first four columns of $Q_x(t)U(t), t=1,\ldots,12$. We observe between four and six dominant peaks; the four clearest ones corresponding to the points 19, 24, 37 and 41 in Figure \ref{fig50points}. This suggests that a reduced dynamical model for the ocean could be built from four to six selected locations. Notice that these four points are on the Pacific ocean, in locations that one would naturally associate with the strongest El Ni\~no signals.

\begin{figure}[h]
\begin{center}
\vspace{-6.0cm}
\includegraphics[height=20cm]{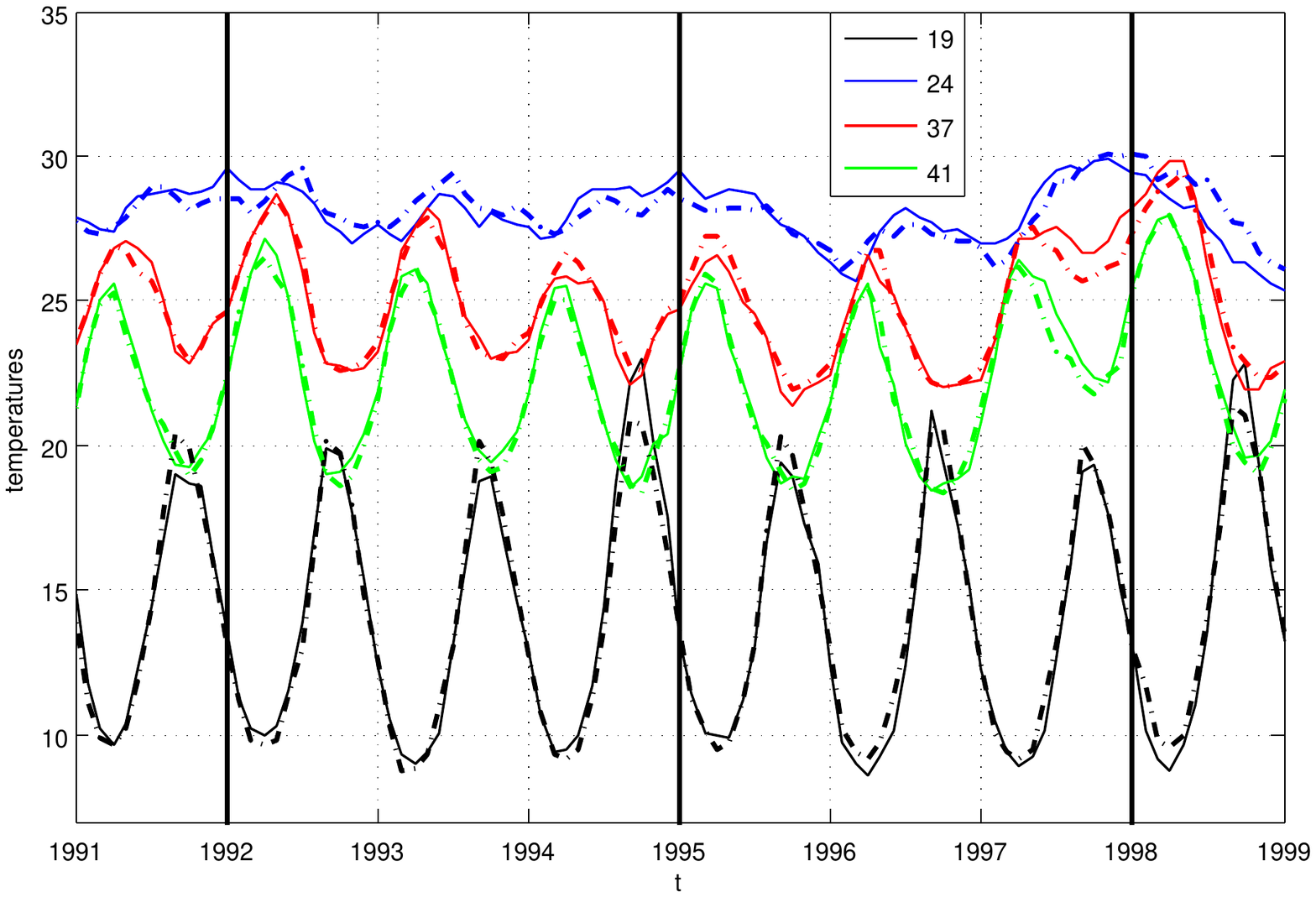}
\vspace{-6.0cm}
\end{center}
\caption{Observed and predicted ocean surface temperature for four selected points. The vertical lines mark El Ni\~no years.}
\label{zvszpred}
\end{figure}

Figure \ref{zvszpred} shows the observed ocean surface temperature for these four points,  comparing them with the ones predicted by the algorithm, in dotted lines. We see that the approximation is quite sharp, particularly near El Ni\~{n}o years, where changes of temperatures are most significant. Even though we have chosen to plot only these four temperatures, all of the 50 points used are well-predicted by the procedure.

\begin{figure}[h]
\begin{center}
\vspace{-6.0cm}
\includegraphics[height=20cm]{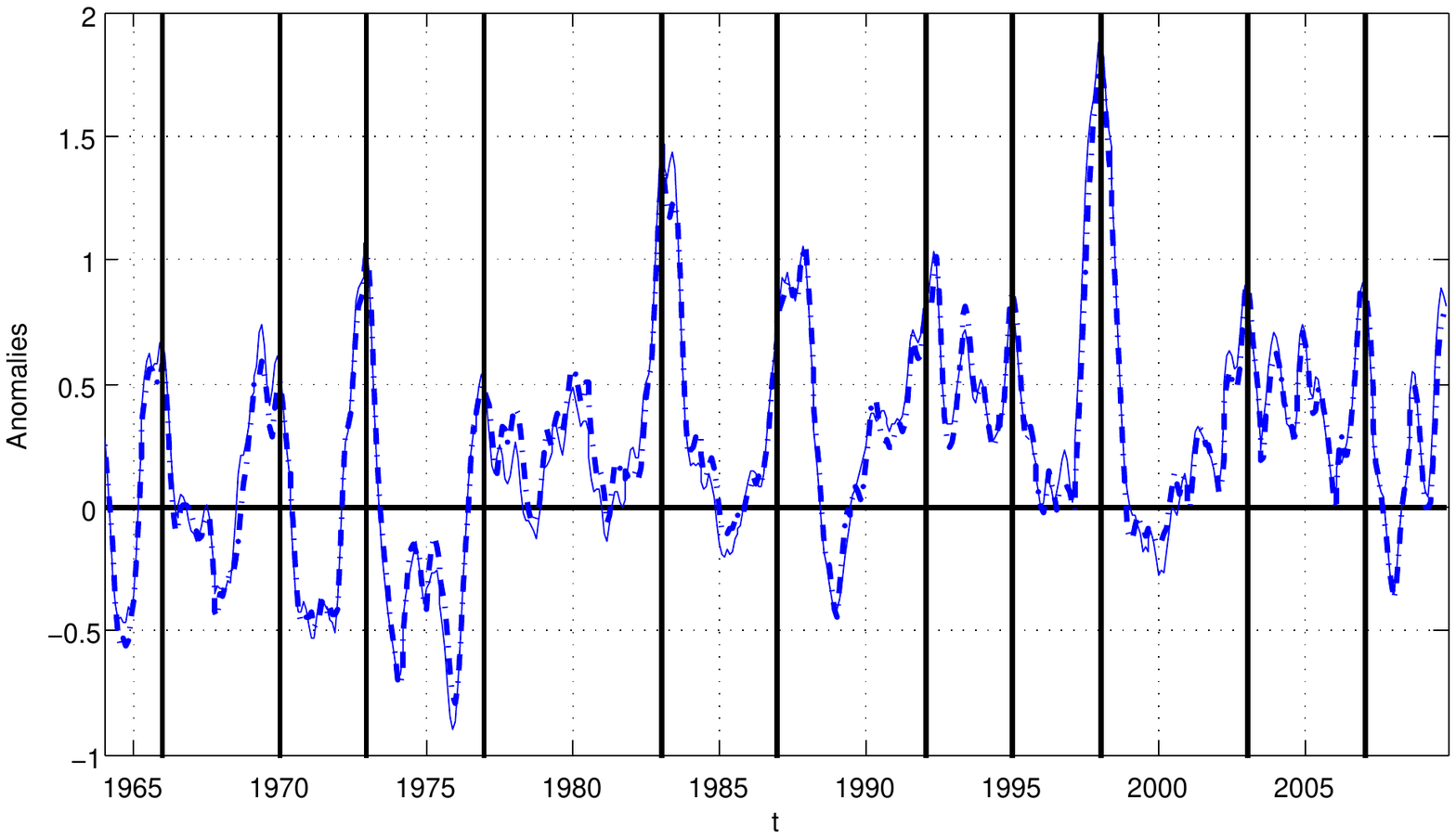}
\vspace{-7.0cm}
\end{center}
\caption{Observed and predicted 3-month mean SST anomalies from February 1964 to October 2009, quantifying El Ni\~no and La Ni\~na intensities. Only El Ni\~no years are marked with vertical lines; La Ni\~na years correspond to strong negative anomalies.}
\label{Anomalies}
\end{figure}

Finally, we monitor the evolution of a measure of the global anomalies associated with El Ni\~no. To this end, we compute a discrete analogue of the running 3-month mean SST anomaly in the El Ni\~{n}o regions \cite{index}. In particular, we average the temperatures on the points 16, 17, 24, 28, 29, 32, 33, 37, 41 and 42 on the map in Figure \ref{fig50points}, for a time window from February 1964 to October 2009. These 10 points are not all strictly included in what are known as El Ni\~{n}o regions (there are four of them, 1+2, 3, 4 and 3.4), but they are the closest on our discrete map to the union of all of them. We observe in Figure \ref{Anomalies} that the warm (positive) peaks coincide with El Ni\~{n}o years. The cold (negative) peaks correspond to La Ni\~{n}a years. In dotted line we have plotted the predicted values of these SST anomalies generated by the principal dynamical component procedure.

\section{Probabilistic perspective and extension to nonlinear dynamics}
\label{Nonlinear}

Throughout this article, we have defined and developed the principal dynamical component procedure in terms of the minimization of a specific cost function: the sum of squares of the prediction errors. In this section, we assign a meaning to this cost in terms of the log-likelihood function of a probabilistic model.
Framing the principal dynamical component procedure in a probabilistic setting has two main advantages: to permit a more thorough interpretation, and to extend its applicability beyond the linear models developed in this article. We sketch such generalization in this section; its algorithmic implementation, under current development, will be presented elsewhere.

Generally, a probabilistic model for a time series $z_j \in R^n$ involves the transition probability density
\begin{equation*}
  T\left(z_{j+1} | z_j \right) .
\end{equation*}
(This corresponds to the Markovian, autonomous scenario, the only one that we address in this section. The extension to non-autonomous and non-Markovian cases, involving a transition probability density of the form $T(z_{j+1} | z_j, z_{j-1}, \ldots, z_{j-r}, t, s)$, is straightforward).  The principal dynamical component proposal considers a dimensional reduction of such transition probability density, using the following elements:
\begin{itemize}
  \item A coordinate system $z = z(x,y)$, $x \in R^m$, $y \in R^{n-m}$, with corresponding projection operators $P_x$ and $P_y$:
  $$x = P_x(z(x,y)), \quad y = P_y(z(x,y)). $$

  \item A reduced dynamical model given by a transition probability density in $R^m$:
  $$ d\left(x_{j+1} | x_j\right). $$

  \item A probabilistic \emph{embedding}
  $$ e(y | x) . $$

  \end{itemize}
The transition probability density for $z$ is then given by
\begin{equation*}
  T\left(z_{j+1} | z_j \right) = J\left(z_{j+1}\right) \, e\left(y_{j+1} | x_{j+1}\right) \, d\left(x_{j+1} | x_{j}\right) ,
\end{equation*}
where $x = P_x(z)$, $y = P_y(z)$, and $J(z)$ is the Jacobian determinant of the coordinate map $z \rightarrow (x,y)$.

A natural measure of the goodness of the model is the log-likelihood function
\begin{equation*}
  L = \sum_{j=1}^{N-1} \log\left[  T\left(z_{j+1} | z_j \right) \right] .
\end{equation*}
In particular, in the setting of Section \ref{Linear}, we have the projections
\begin{equation}
 P_x(z) = Q_x' z , \quad P_y(z) = Q_y' z ,
 \label{projections}
\end{equation}
where $Q = [Q_x Q_y]$ is orthogonal, so $J(z) = 1$. The embedding and reduced dynamics are given by the isotropic Gaussians
\begin{equation}
 e(y|x) = \mathcal{N}(0, \sigma^2 I_{N-m})
 \label{embiso}
\end{equation}
and
\begin{equation}
 d\left(x_{j+1} | x_j\right) =  \mathcal{N}(A x_j, \sigma^2 I_m) ,
 \label{dyniso}
\end{equation}
where $I_k$ stands for the $k\times k$ identity matrix. Consequently, the log-likelihood function is given by
\begin{equation*}
  L = \sum_{j=1}^{N-1}  -\left[\frac{n}{2} \log(2\pi)+ n \log(\sigma) + \frac{1}{2\sigma^2} \left(\| x_{j+1} - A x_j \|^2 + \|y_{j+1}\|^2 \right) \right] .
  \label{LLiso}
\end{equation*}
Thus maximizing the log-likelihood $L$ over $Q$ and $A$ is equivalent to minimizing the cost function
\begin{equation*}
  c = \frac{1}{N-1} \sum_{j=1}^{N-1}  \left(\| x_{j+1} - A x_j \|^2 + \|y_{j+1}\|^2 \right)\, ,
\end{equation*}
that we have used throughout the paper; the corresponding optimal value of $\sigma$ is given by
\begin{equation*}
  \sigma = \left( \frac{c}{n} \right)^{\frac{1}{2}} .
\end{equation*}

This interpretation immediately suggests the following generalization, which remains within the realm of Gaussian distributions and linear maps: keep the orthogonal projections in (\ref{projections}), but replace the embedding (\ref{embiso}) and dynamical model (\ref{dyniso}) by the more general
\begin{eqnarray*}
 e(y|x) &=& \mathcal{N}(0, \Sigma_y) , \\
 d\left(x_{j+1} | x_j\right) &=&  \mathcal{N}(A x_j, \Sigma_x) ,
\end{eqnarray*}
where $\Sigma_x$ and $\Sigma_y$ are general covariance matrices. The resulting log-likelihood function is
\begin{equation*}
  L = \sum_{j=1}^{N-1}  -\frac{1}{2}\left[ \log((2\pi)^n |\Sigma_x| |\Sigma_y|) + \left( x_{j+1} - A x_j, \Sigma_x^{-1} (x_{j+1} - A x_j) \right) + \left(y_{j+1},\Sigma_y^{-1} y_{j+1}\right) \right] .
  \label{LLnoniso}
\end{equation*}
This formulation has the advantage of providing a natural ranking of the coordinates $x$ and $y$, through the principal components of the corresponding covariance matrices.

More generally, one can propose different, typically nonlinear, families of distributions, projections and dynamical models, and maximize the corresponding log-likelihood function. The proposed distributions can be given parametrically, in which case the maximization of the log-likelihood is over their parameters, or non-parametrically, for instance as an extension of the methodology proposed in \cite{DE}. Thus the principal dynamical component methodology extends naturally to very general scenarios, with nonlinear reduced dynamical manifolds, stochastic, nonlinear dynamical models, and non-Gaussian embeddings. This extension, however, goes beyond the scope of this paper, and will be pursued elsewhere.

\section{Conclusions}
\label{Conclusion}

A new methodology has been developed for the dimensional reduction of time series. The procedure seeks a low dimensional manifold $x$ and a dynamical model $x_{j+1} = D(x_j, x_{j-1}, \ldots, t)$ that minimize the predictive uncertainty of the series.  The procedure has been  successfully tested on synthetic data, and illustrated with a real application to time series of sea-surface temperature over the ocean. Finally, a probabilistic interpretation of  the principal  dynamical component procedure was proposed, providing a conceptual extension to general nonlinear, non-Gaussian settings.

\section{Acknowledgements}
The work of M. D. de la Iglesia is  partially supported by D.G.E.S, ref. BFM2006-13000-C03-01, Junta de Andaluc\'{i}a, grants FQM-229, FQM-481, P06-FQM-01738 and Subprograma de estancias de movilidad posdoctoral en el extranjero, MICINN, ref. -2008-0207, and that of E. G. Tabak is partially supported by the National Science Foundation under grant number DMS 0908077.

\bibliographystyle{elsart-num}

\begin{thebibliography}{99}


\bibitem{Jit} Jolliffe, I. T.,  \emph{Principal component analysis}, Springer series in Statistics, Springer-Verlag,
1986.

\bibitem{EOF} Lorenz, E. N., ``Empirical orthogonal functions and statistical weather prediction'', Statistical forecast project report 1, Dept. of Meteor., MIT, 1956.


\bibitem{DE} Tabak, E. and Vanden-Eijnden, E.,
``Density estimation by dual ascent of the log-likelihood'',
\emph{Comm. Math. Sci.}, {\bf 8} , 217-233, 2010.

\bibitem{IRI} The International Research Institute for Climate and Society, http://iridl.ldeo.columbia.edu/SOURCES/.NOAA/.NCDC/.ERSST/

\bibitem{index} Trenberth, K. E., and D. P. Stepaniak, ``Indices of El Ni\~no evolution'', \emph{ J. Climate}, {\bf 14}, 1697-1701, 2001.

\bibitem{Wei} Wei, W. W.,  \emph{Time series. Univariate and multivariate methods}, Addison-Wesley Publishing Company, Advanced Book Program, Redwood City, CA, 1990.



\end{thebibliography}

\end{document}